\documentstyle[11pt]{article}
\input epsf.tex
\font\goth=eufm10

\newtheorem{theorem}{Theorem}

\newtheorem{proposition}[theorem]{Proposition}
\newtheorem{lemma}[theorem]{Lemma}
\newtheorem{example}[theorem]{Example}

\def\boxit#1#2{\setbox1=\hbox{\kern#1{#2}\kern#1}%
\dimen1=\ht1 \advance\dimen1 by #1 \dimen2=\dp1 \advance\dimen2 by #1
\setbox1=\hbox{\vrule height\dimen1 depth\dimen2\box1\vrule}%
\setbox1=\vbox{\hrule\box1\hrule}%
\advance\dimen1 by .4pt \ht1=\dimen1
\advance\dimen2 by .4pt \dp1=\dimen2 \box1\relax}

\def\Proof{\medskip\noindent {\it Proof --- \ }}

\def\cqfd{\hfill $\Box$ \bigskip}
\def\adots{\mathinner{\mkern2mu\raise1pt\hbox{.}
\mkern3mu\raise4pt\hbox{.}\mkern1mu\raise7pt\hbox{.}}}
\def\<{\lngle\,}
\def\>{\,\rangle}

\def\ie{{\it i.e. }}
\def\eg{{\it e.g. }}

\def\SG{\hbox{\goth S}}

\def\N{{\bf N}}
\def\Z{{\bf Z}}

\def\T2{{\cal T}}
\def\F{{\cal F}}
\def\G{{\cal G}}
\def\A{{\cal A}}

\def\SS{{\cal S}}

\def\slchap{\widehat{sl}}

\def\ASG{\widetilde{\hbox{\goth S}}}
\def\EASG{\widehat{\hbox{\goth S}}}
\def\AH{\widetilde{H}}
\def\EAH{\widehat{H}}

\newdimen\Squaresize \Squaresize=14pt
\newdimen\Thickness \Thickness=0.5pt

\def\Square#1{\hbox{\vrule width \Thickness
   \vbox to \Squaresize{\hrule height \Thickness\vss
      \hbox to \Squaresize{\hss#1\hss}
   \vss\hrule height\Thickness}
\unskip\vrule width \Thickness}
\kern-\Thickness}
 
\def\Vsquare#1{\vbox{\Square{$#1$}}\kern-\Thickness}

\def\today{\number\day \space\ifcase\month\or
 Janvier\or F\'evrier\or Mars\or Avril\or Mai\or Juin\or
 Juillet\or Ao\^ut\or Septembre\or Octobre\or Novembre\or D\'ecembre\fi
 \space\number\year}

\title{\bf Decomposition numbers and canonical bases} 
\author{\rm Bernard {\sc Leclerc} \\ [2mm]
\footnotesize \it D\'epartement de Math\'ematiques,\\
\footnotesize \it Universit\'e de Caen, 14032 Caen cedex, France.\\
\footnotesize\tt leclerc@math.unicaen.fr \\
\footnotesize\tt \qquad http://matin.math.unicaen.fr/\~{}leclerc}
\date{}

\begin{document}
\maketitle

\vskip 0.5cm

\begin{abstract}\noindent
We obtain some simple relations between decomposition
numbers of quantized Schur algebras at an $n$th root
of unity (over a field of characteristic 0).
These relations imply that every decomposition number for
such an algebra occurs as a decomposition number for some Hecke
algebra of type $A$.
We prove similar relations between coefficients of the 
canonical basis of the 
$q$-deformed Fock space representation of $U_q(\slchap_n)$
introduced in \cite{LT}. 
It follows that these coefficients can all be expressed in terms 
of those of the global crystal basis of the irreducible 
sub-repre\-sen\-ta\-tion generated by the vacuum vector.

As a consequence, using works of Ariki and Varagnolo-Vasserot,
it is possible to give a new proof of Lusztig's character formula for the
simple $U_v(sl_r)$-modules at roots of unity, which does not 
involve representations of $\slchap_r$ of negative level.
\end{abstract}
 

\section{Introduction}
\label{SECT1}
It is well known that the $p$-modular decomposition numbers
for the symmetric groups are equal to certain composition multiplicities
of the Weyl modules of the general linear groups
over an algebraically closed field of characteristic $p$
(see \eg \cite{JK}, p.341, \cite{Ja}, p.317).
This follows from the properties of the Schur functor
from representations of the general linear groups
to representations of the symmetric groups.
In other words, every decomposition number for a symmetric
group $\SG_m$ occurs as a decomposition number for the 
associated Schur algebra $\SS_m$.

More recently Erdmann \cite{Erd1} has obtained a kind of converse
result, namely that the decomposition matrix of $\SS_m$
is a sub-matrix of the decomposition matrix of the symmetric
group $\SG_{m'}$ of rank $m'= pm+(p-1){m\choose 2}$.
The proof uses the Frobenius map of the general linear
group and the theory of tilting modules.

The first of these results has a quantum analogue
given by James \cite{Ja1} and Dipper-James~\cite{DJ}. 
It states that the decomposition matrix of the Hecke  
algebra $H_m(v)$ of type $A_{m-1}$ embeds into the decomposition 
matrix of the $v$-Schur algebra $\SS_m(v)$.
But, as noted in \cite{DPS}, ``there is no evident quantization''
of Erdmann's theorem, because the quantum analogue of the
Frobenius map does not go from the quantum algebra $U_v(gl_r)$ to
itself but from $U_v(gl_r)$ to
the classical enveloping algebra $U(gl_r)$
(see also \cite{Do}, p. 104).

Nevertheless, we will show in the first part of this note that conversely
every decomposition number for $\SS_m(v)$
occurs as a decomposition number for the Hecke algebra 
$H_{m''}(v)$, where $m'' = m + 2(n-1){m\choose 2}$ and
$n$ is the multiplicative order of $v^2$. 
Our formula also relies on some properties of the tilting modules,
but the argument of Erdmann based on the Frobenius map is replaced
by a result of Andersen relating  
injective $U_v(sl_r)$-modules to tilting modules
(\cite{An}, Prop. 5.8).
This result being valid only over a field of characteristic 0,
the same restriction applies to our formula.
So the question remains open in the so-called mixed case.

In the second part of this note we consider the coefficients
of the canonical basis introduced in \cite{LT} of the Fock space 
representation $\F$ of $U_q(\slchap_n)$.
These are polynomials $d_{\lambda,\mu}(q)$ indexed by pairs $(\lambda,\mu)$
of partitions, as the decomposition numbers $d_{\lambda,\mu}$
of the $v$-Schur algebras.
The main conjecture of \cite{LT} was that 
\begin{equation}\label{CONJ}
d_{\lambda,\mu}(1) = d_{\lambda,\mu}\,, 
\end{equation}
\ie the $d_{\lambda,\mu}(q)$ are
some $q$-analogues of the decomposition numbers.
In \cite{VV}, Varagnolo and Vasserot proved that the
$d_{\lambda,\mu}(q)$ are certain parabolic Kazhdan-Lusztig polynomials
for the affine symmetric groups $\ASG_r$, and thus showed 
that Eq.~(\ref{CONJ}) results from the so-called Lusztig conjecture
on the characters of the simple $U_v(sl_r)$-modules \cite{Lu5}.
This last conjecture is proved, due to work of Kazhdan-Lusztig \cite{KL}
and Kashiwara-Tanisaki \cite{KT} involving 
a certain category of $\slchap_r$-modules of negative level $-n-r$.

In Section~\ref{SECT4} we prove a relation between the $d_{\lambda,\mu}(q)$
which is a $q$-analogue of the relation between decomposition
numbers established in the first part. 
Its proof is based on a theorem of Soergel on Kazhdan-Lusztig
polynomials (\cite{So1}, Theorem 5.1)
and on a certain duality satisfied by the canonical basis of $\F$,
which was announced in \cite{LT} and proved in \cite{LT2} (Theorem 7.14).

This implies that all the coefficients in the expansion 
of the canonical basis of $\F$
on its natural basis consisting of the pure $q$-wedges
can be expressed in terms of those of the global crystal basis of
the simple submodule $L(\Lambda_0)$ generated by the vacuum vector.

\bigskip
It follows from these two sets of relations
(between decomposition numbers on one hand, and between 
coefficients of the canonical basis on the other hand)
that it is enough for proving Eq.~(\ref{CONJ})
to verify it in the case where $\mu$ is an $n$-regular partition.
This was done by Ariki in \cite{Ar}, using results of
Kazhdan-Lusztig and Ginzburg on representations of
affine Hecke algebras, and the geometric description by Lusztig
of the canonical basis of $U_q^-(\slchap_n)$.
Therefore, as shown in Section~\ref{SECT5}, Ariki's theorem implies Lusztig's formula 
for the characters of the simple $U_v(sl_r)$-modules, as well
as Soergel's formula for the characters of the indecomposable
tilting $U_v(sl_r)$-modules \cite{So2}. 

Varagnolo and Vasserot have already attempted in \cite{VV} to give
a proof of Lusztig's formula based on the level $1$ Fock space representation $\F$
of $U_q(\slchap_n)$. 
They have extended the action of $U_q^-(\slchap_n)$ on $\F$
to an action of the Hall algebra of the cyclic quiver, and have
defined in this way, using intersection cohomology, a geometric
basis ${\bf B}$ of $\F$. 
They conjectured that ${\bf B}$ coincides with the canonical basis
of \cite{LT}, which is defined in an algebraic way, and 
proved that this equality would amount to a $q$-analogue of 
Lusztig's formula. 
This very interesting conjecture is proved by Schiffmann in
the case $n=2$ \cite{OS}, but remains open in general.


\section{Statement of results}
\label{SECT2}

We fix two integers $m,n \ge 2$ and denote by $D_m$ the decomposition matrix of
the quantized Schur algebra $\SS_m(v)$ over a field of characteristic
0, with parameter $v$ such that $v^2$ is a primitive $n$th root of unity.
(There are several ways of introducing the $v$-Schur algebra. 
Here we regard $\SS_m(v)$ as a suitable quotient of the 
quantized enveloping
algebra $U_v(gl_r)$ for $r\ge m$, as in \cite{Du}.)
The rows and columns of $D_m$ are labelled by the set of partitions
of $m$, and we use the notational convention of \cite{Jam}, that is,
the entry on row $\lambda$ and column $\mu$ is
\begin{equation}\label{DEFDC}
d_{\lambda,\mu} := [W(\lambda') : L(\mu')]_{\SS_m(v)}\,,
\end{equation}
where $\lambda'$ stands for the partition conjugate to $\lambda$,
and $W(\lambda)$, $L(\lambda)$ are respectively the Weyl module 
and the simple module with highest weight $\lambda$.

We fix some $r\ge 2$ and put $\rho_r = (r-1,r-2,\ldots ,1,0)$.
A partition $\lambda$ with at most $r$ parts is identified with an $r$-tuple
in $\N^r$ in the standard way by appending a tail of $0$.
Given two such partitions $\lambda,\,\mu$ of $m$, we shall define
two partitions $\widehat\mu$ and $\widetilde\lambda$
of $m''=m + (n-1)r(r-1)$ in the following way.
There is a unique decomposition $\mu = \mu^{(0)} + n\mu^{(1)}$
with the partition $\mu^{(0)}$ $n$-restricted, that is,
the difference between any two consecutive parts is $< n$.
Define
\begin{eqnarray}
\widehat\mu &=& 2(n-1)\rho_r + w_0(\mu^{(0)}) + n\mu^{(1)}\,,
\label{DEFMUCHAP}\\
\widetilde\lambda &=& \lambda + ((n-1)(r-1),\ldots ,(n-1)(r-1))\,, 
\end{eqnarray}
where for $\alpha = (\alpha_1,\alpha_2,\ldots ,\alpha_r)$ we set
$w_0(\alpha) =  (\alpha_r,\alpha_{r-1},\ldots ,\alpha_1)$.
It is easy to check that the parts of $\widehat\mu$ are pairwise
distinct, and thus $\widehat\mu$ is always an $n$-regular
partition, \ie no part has multiplicity $\ge n$.

\begin{theorem}\label{TH1}
For all partitions $\lambda,\,\mu$ of $m$ of length $\le r$, we have
\begin{equation}
d_{\lambda',\,\mu'}=d_{\widetilde\lambda,\widehat\mu}\,.
\end{equation}
\end{theorem}
Since $\widehat\mu$ is $n$-regular, we see that 
$d_{\widetilde\lambda,\widehat\mu}$ is a decomposition number
for the Hecke algebra $H_{m''}(v)$, and therefore if $r\ge m$ every
decomposition number for $\SS_m(v)$ occurs as a decomposition
number for $H_{m''}(v)$.

Let now $d_{\lambda,\mu}(q)$ be the polynomial 
defined in \cite{LT} (see also \cite{LT2}, Section 7).
The affine symmetric group $\ASG_r$ acts on $\Z^r$
via its level $n$ action $\pi_n$ (see \cite{LT2}, Section 2),
and the $q$-analogue of the linkage principle states that
$d_{\lambda,\mu}(q) \not = 0$ only if $\lambda+\rho_r$ and $\mu+\rho_r$ belong
to the same orbit under $\ASG_r$, \ie the partitions $\lambda$
and $\mu$ have the same $n$-core.
Let $\nu$ be the unique point of the orbit $\ASG_r(\mu+\rho_r)$
such that $\nu_1\ge\nu_2\ge \cdots \ge\nu_r$ and
$\nu_1-\nu_r\le n$.
The stabilizer of $\nu$ is a standard parabolic subgroup 
of $\ASG_r$ and we write $\ell_\mu$ for the length of its
longest element. 
We also put $\ell = r(r-1)/2$, the length of the longest element
of $\SG_r$.

\begin{theorem}\label{TH2}
For all partitions $\lambda,\,\mu$ of $m$ of length $\le r$, we have
\begin{equation} 
d_{\lambda',\,\mu'}(q)=q^{\ell-\ell_\mu}\,
d_{\widetilde\lambda,\widehat\mu}(q^{-1})\,.
\end{equation} 
\end{theorem}
Since $\widehat\mu$ is $n$-regular, it labels
a vector of the global crystal basis of the simple module $L(\Lambda_0)$. 
Hence every vector 
\begin{equation}
{\cal G}^+_{\mu'}=\sum_\lambda d_{\lambda',\,\mu'}(q)\, |\lambda'\rangle
\end{equation}
of the canonical basis of the Fock
space $\F$ can be easily computed from the corresponding vector 
${\cal G}^+_{\widehat\mu}$ of the global basis of $L(\Lambda_0)$.

\begin{example}{\rm

Let $n=3$ and $\mu=(6,2,1)$, a partition of $m=9$.
The non-zero $d_{\lambda',\,\mu'}$ and $d_{\lambda',\,\mu'}(q)$
are obtained for 
$
\lambda = (6,2,1),\ (7,1,1),\ (6,3),\ (8,1).
$
We can take $r=3$, so that $m''=9+2.3.2 = 21$.
We have
\begin{eqnarray*}
\mu &=& \mu^{(0)} + 3\mu^{(1)} = (3,2,1) +3(1,0,0)\,,\\
\widehat{\mu} &=& (1,2,3)+4(2,1,0)+3(1,0,0) = (12,6,3)\,.
\end{eqnarray*}
Hence, taking
$ 
\widetilde{\lambda} = (10,6,5),\ (11,5,5),\ (10,7,4),\ (12,5,4) , 
$
respectively, we have
\[
d_{\lambda',\,(3,2,1,1,1,1)} = d_{\widetilde{\lambda},(12,6,3)}\,.
\]
Moreover, $\mu+\rho_r = (8,3,1)$ is regular for the level 3 action
of $\ASG_3$, hence $\ell_\mu = 0$ and we have
\[ 
d_{\lambda',\,(3,2,1,1,1,1)}(q) = q^3\,d_{\widetilde{\lambda},(12,6,3)}(q^{-1})\,. 
\]
}
\end{example}


\section{Proof of Theorem~1}
\label{SECT3}

Let $m, r$ be two integers $\ge 2$.
Let $U_v(gl_r)$ and $U_v(sl_r)$ denote Lusztig's restricted
specialization at $v$ of the quantum enveloping algebras
of $gl_r$ and $sl_r$ respectively.
As mentioned above, we regard the $v$-Schur algebra $\SS_{m,r}(v)$
as the homomorphic image of $U_v(gl_r)$
obtained via its action on the $m$th tensor power of the vector 
representation, and we put $\SS_m(v) = \SS_{m,m}(v)$
(see \cite{Du}, \cite{DPS}).
Thus every $\SS_{m,r}(v)$-module becomes a $U_v(gl_r)$-module
and by restriction a $U_v(sl_r)$-module. 
We shall abuse notation and denote in the same way these various modules.
The simple modules and Weyl modules of $\SS_{m,r}(v)$ are
labelled by the set of partitions of $m$ with at most $r$ parts,
and are denoted by $L(\lambda)$ and $W(\lambda)$.
It is known that for such partitions $\lambda,\, \mu$,
\begin{equation}
[W(\lambda) : L(\mu)]_{\SS_{m,r}(v)} = 
[W(\lambda) : L(\mu)]_{\SS_m(v)}\,.
\end{equation}
Hence, it is enough to determine the decomposition numbers 
$d_{\lambda,\mu}$ as defined in (\ref{DEFDC}).

Let now $\lambda,\,\mu$ be two partitions of $m$ and fix an integer 
$r$ greater or equal to the number of parts of both $\lambda$ and
$\mu$, for example $r = m$. We have
\begin{equation}
d_{\lambda',\,\mu'} 
= [W(\lambda) : L(\mu)]_{\SS_{m,r}(v)}
=[W(\lambda) : L(\mu)]_{U_v(gl_r)}
=[W(\lambda) : L(\mu)]_{U_v(sl_r)}\,.
\end{equation}
Let $I(\mu)$ denote the injective hull of the $U_v(sl_r)$-module
$L(\mu)$. 
By the reciprocity formula (see \cite{APW2}) we have
\begin{equation}
[W(\lambda) : L(\mu)]_{U_v(sl_r)}
= [I(\mu) : W(\lambda)]_{U_v(sl_r)} \,.
\end{equation}
Given a partition $\alpha$ of $m$ of length $\le r$,
let $T(\alpha)$ denote the unique indecomposable tilting $\SS_{m,r}(v)$-module
such that $\alpha$ is the maximal partition $\beta$ (for the usual dominance order)
for which $[T(\alpha) : L(\beta)] \not = 0$.
Consider the $\SS_{m'',r}(v)$-module $T(\widehat\mu)$, where 
$\widehat\mu$ is defined by (\ref{DEFMUCHAP}).
Then, regarding $T(\widehat\mu)$ as a ${U_v(sl_r)}$-module,  
we have by \cite{An}, Prop. 5.8 that $T(\widehat\mu)$ is isomorphic to $I(\mu)$.
(In \cite{An} there are some restrictions on $n$, that is,
on the multiplicative order of the root of unity $v^2$.
These restrictions have been later removed in \cite{AnP},
and although the statement of \cite{An}, Prop. 5.8 is not given
in that paper, the arguments easily carry over.)
Thus, noting that 
the $\SS_{m,r}(v)$-module $W(\lambda)$ and the $\SS_{m'',r}(v)$-module
$W(\widetilde\lambda)$ become isomorphic when considered as $U_v(sl_r)$-modules,
we get
\begin{equation}
d_{\lambda',\,\mu'} = [T(\widehat\mu) : W(\widetilde\lambda)]_{U_v(sl_r)} 
= [T(\widehat\mu) : W(\widetilde\lambda)]_{\SS_{m'',r}(v)}
= [T(\widehat\mu) : W(\widetilde\lambda)]_{\SS_{m''}(v)}\,.
\end{equation}
Finally, by \cite{DPS}, Proposition 8.2 (a),
or \cite{Do}, Proposition 4.15 (ii), 
\begin{equation}
[T(\widehat\mu) : W(\widetilde\lambda)]_{\SS_{m''}(v)}
= [ W(\widetilde\lambda') : L(\widehat\mu')]_{\SS_{m''}(v)}
= d_{\widetilde\lambda,\widehat\mu}\,,
\end{equation}
and Theorem~\ref{TH1} is proved.


\section{Proof of Theorem~2}
\label{SECT4}
We shall use the same notation as in \cite{LT2}.
In particular $\SG_r$, $\ASG_r$ and $\EASG_r$ stand for the symmetric group, 
the affine symmetric
group and the extended affine symmetric group, respectively, and 
$H_r$, $\AH_r$ and $\EAH_r$ are the corresponding Hecke algebras.
The standard generators of $\EASG_r$ are denoted by $s_0,\ldots , s_{r-1}, \tau$,
and $\SG_r$, $\ASG_r$ are the subgroups generated by $s_1,\ldots , s_{r-1}$
and $s_0,\ldots , s_{r-1}$, respectively.
We write $w_0$ for the longest element of~$\SG_r$.
For $x,w \in \EASG_r$, we have the Kazhdan-Lusztig polynomial $P_{x,w}(q)$,
and the inverse Kazhdan-Lusztig polynomial $Q_{x,w}(q)$ defined
via the equations
\begin{equation}
\sum_{x\in \EASG_r} Q_{x,z}(-q)P_{x,w}(q) = \delta_{z,w}\,,
\qquad (z,\,w\in \EASG_r)\,.
\end{equation}
The group $\EASG_r$ acts on $P_r := \Z^r$ via the level $-n$ action $\pi_{-n}$, 
giving rise to the parabolic Kazhdan-Lusztig polynomials
$P^-_{\mu,\lambda}(q)$ indexed by $\lambda, \mu \in P_r$.
Let $Q^-_{\mu,\lambda}(q)$ denote the inverse parabolic Kazhdan-Lusztig
polynomial defined by the equations
\begin{equation}
\sum_{\mu\in P_r} Q^-_{\mu,\alpha}(-q)P^-_{\mu,\lambda}(q) = \delta_{\alpha,\lambda}\,
\qquad (\alpha,\,\lambda \in P_r)\,.
\end{equation}
We recall (see \cite{LT2}, Section 2) that if $\lambda_i > \lambda_{i+1}$ then 
\begin{equation}
P^-_{s_i\mu,\lambda}(q) = \left\{
\matrix{0 & \mbox{if $\mu_i = \mu_{i+1}$,}\cr
qP^-_{\mu,\lambda}(q) & \mbox{if $\mu_i > \mu_{i+1}$.}
}
\right.
\end{equation}
Thus if $\lambda$ is strictly dominant, \ie
$\lambda \in P^{++}_r := 
\{\alpha \in P_r \, | \, \alpha_1 > \alpha_2 > \cdots > \alpha_r \}$,
we get that $P^-_{\mu,\lambda}(q) = 0$ for $\mu \not \in \SG_rP^{++}_r$
(\ie if the coordinates of $\mu$ are not pairwise distinct),
and otherwise
\begin{equation}\label{REFL}
P^-_{s\beta , \lambda}(q) = q^{\ell(s)}\,P^-_{\beta,\lambda}(q)\,,
\qquad (\beta \in P^{++}_r, s\in \SG_r)\,.
\end{equation} 

Let us fix $\lambda,\,\alpha \in P^{++}_r$. Using (\ref{REFL}) we obtain 
\begin{eqnarray*}
\delta_{\alpha,\lambda} &=& \sum_{\mu\in P_r} Q^-_{\mu,\alpha}(-q)P^-_{\mu,\lambda}(q) \\
&=&\sum_{\beta \in P^{++}_r}\left(
\sum_{s\in\SG_r} Q^-_{s\beta,\alpha}(-q)P^-_{s\beta,\lambda}(q) \right) \\
&=&\sum_{\beta \in P^{++}_r}\left(\sum_{s\in\SG_r} q^{\ell(s)}
Q^-_{s\beta,\alpha}(-q)\right)P^-_{\beta,\lambda}(q)\,.
\end{eqnarray*}
Hence, setting 
\begin{equation}\label{DEFr}
r_{\beta,\alpha}(q) := \sum_{s\in\SG_r} (-q)^{\ell(s)}Q^-_{s\beta,\alpha}(q) \,,
\qquad (\alpha,\,\beta \in P^{++}_r)\,,
\end{equation}
we have obtained
\begin{lemma} For $\alpha,\,\lambda \in P^{++}_r$, there holds
\begin{equation}\label{PARABINV}
\sum_{\beta \in P^{++}_r} r_{\beta,\alpha}(-q)P^-_{\beta,\lambda}(q) = 
\delta_{\alpha,\lambda}\,.
\end{equation}
\cqfd
\end{lemma}
Now comparing Eq.~(\ref{PARABINV}) with \cite{LT2}, Eq. (91) and Corollary 7.15
we get 
\begin{proposition}\label{P1}
Let $\lambda,\,\mu$ be two partitions of $m$ with at most 
$r$ parts, and put $\beta = \lambda+\rho_r$, $\alpha=\mu+\rho_r$.  Then
\[
r_{\beta,\alpha}(q) = d_{\lambda',\mu'}(q) \,.
\]
\cqfd
\end{proposition}

Let $W^f$ denote the set of minimal length representatives for
the right cosets $\SG_r \backslash \ASG_r$. 
We shall express the polynomials $r_{\beta,\alpha}(q)$
in terms of the 
\begin{equation}
m^{x,w}(q) := \sum_{s\in\SG_r} 
(-q)^{\ell(w_0) - \ell(s)}\,
Q_{sx,w_0w}(q)\,,\qquad (x,w \in W^f)\,.
\end{equation}
These are the inverse parabolic Kazhdan-Lusztig polynomials 
coming from the right action of $\AH_r$ on the parabolic
module ${\bf 1}_q \otimes_{H_r} \AH_r$, where
${\bf 1}_q$ denotes the 1-dimensional right $H_r$-module
given by ${\bf 1}_q T_i = q{\bf 1}_q\,, \ (1\le i \le r-1)$,
(see \cite{So1}, Proposition 3.7).

We need to recall some notation from \cite{LT2}.
Let $k\in \Z^*$. We write $\A_{r,k}$ for the fundamental alcove
in $P_r$ associated with the level $k$ action $\pi_k$ of $\EASG_r$.
For $\lambda\in P_r$, we denote by $w(\lambda,k)$ the element of 
$\EASG_r$ of minimal length such that $w(\lambda,k)^{-1} \lambda \in \A_{r,k}$.

In \cite{So1}, Soergel works with $\ASG_r$ and $\AH_r$ instead of 
$\EASG_r$ and $\EAH_r$.
Recall that every element $w\in\EASG_r$ can be expressed uniquely as
$w=\sigma\tau^a,\ (\sigma\in\ASG_r,\,a\in\Z)$.
We define a projection of $\EASG_r$ onto $\ASG_r$ by setting
$\underline{w} := \sigma$. 

Let $\mu = \sigma\lambda\in P_r$ with $\sigma\in \ASG_r$. Then clearly
$\sum_{i=1}^r \mu_i = \sum_{i=1}^r \lambda_i$, hence $\ASG_r$ acts
on the hyperplane $\sum_{i=1}^r \lambda_i =0$ that we may identify
with $\underline{P}_r := P_r/\Z(1,\ldots,1)$.
We denote by $\lambda \mapsto \underline{\lambda}$ the 
natural projection $P_r \rightarrow \underline{P}_r$,
and we take 
\[
\underline{\A}_{r,n} := \{\underline{\lambda} 
 \ |\  \lambda_1 \ge \cdots \ge \lambda_r,\ \lambda_1-\lambda_r \le n\}
\]
as fundamental alcove of $\ASG_r$ acting on $\underline{P}_r$ via the
action $\underline{\pi}_n$ induced by $\pi_n$.
Let $\underline{w}(\lambda,n)$ denote the element of $\ASG_r$ of
minimal length such that 
$\underline{w}(\lambda,n)^{-1} \underline{\lambda}\in\underline{\A}_{r,n}$.
Then one has $\underline{w}(\lambda,n) = \underline{w(\lambda,n)}$.
We shall also use the shorthand notation $w_\alpha := \underline{w}(\alpha,n)$.

\begin{proposition}\label{PROP5}
For $\alpha,\,\beta \in P^{++}_r$, there holds
\[
r_{\beta,\alpha}(q) =  m^{w_\beta , w_\alpha}(q) \,.
\]
\end{proposition}
\Proof
First we note that if $\alpha\in P^{++}_r$ then $w_\alpha \in W^f$.
Indeed, since the alcove $A=w_\alpha\,\underline{\A}_{r,n}$
contains $\underline{\alpha}\in \underline{P}^{++}_r$, it has to be
contained in $\underline{P}^{+}_r$, which is equivalent to
$w_\alpha \in W^f$.

Let us rewrite the definition (\ref{DEFr}) of $r_{\beta,\alpha}$
in terms of the (ordinary) inverse Kazhdan-Lusztig polynomials
$Q_{x,w}$.
By \cite{So1}, Proposition 3.7, we have
\[
Q^-_{\mu,\lambda} = Q_{w(\mu,-n)\,,\,w(\lambda,-n)}\,,
\qquad (\lambda , \mu \in P_r)\,.
\]
We have to relate the level $-n$ and $+n$ actions of $\EASG_r$
and $\ASG_r$ on $P_r$.
Let $\sharp$ denote the automorphism of $\EASG_r$ defined by
$\tau^\sharp = \tau^{-1}$ and $s_i^\sharp = s_{-i}$, where $i$ is
understood modulo $r$.
Then it is easy to check that 
\[
\underline{w(s\beta,-n)} = (sw_0).\underline{w(\beta,n)}^\sharp
= sw_0\,.\,w_\beta^\sharp\,,
\qquad (\beta \in P^{++}_r,\ s \in \SG_r)\,.
\]
Therefore, noting that $\SG_r$ is stable under $\sharp$
and $Q_{x^\sharp,y^\sharp} = Q_{x,y}$, we
get for $\alpha,\,\beta \in P^{++}_r$, 
\[
r_{\beta,\alpha}(q) =  \sum_{s\in\SG_r} (-q)^{\ell(s)}
Q_{sw_0w_\beta\,,\,w_0w_\alpha}(q) 
=
\sum_{s'\in\SG_r} (-q)^{\ell(w_0)-\ell(s')}
Q_{s'w_\beta\,,\,w_0w_\alpha}(q)
=m^{w_\beta , w_\alpha}(q)\,.
\]

\cqfd

In \cite{So1}, Soergel considers a map $A \mapsto \widehat A$ on the 
set of $\ASG_r$-alcoves. The next lemma relates this operation to the
map $\mu \mapsto \widehat\mu$ on partitions defined in Section~\ref{SECT2}.
\begin{lemma}\label{LEMM}
Let $\mu$ be a partition with at most $r$ parts, and put 
$\alpha := \mu + \rho_r \in P^{++}_r$.
Let $\xi$ be the point in $\underline{\A}_{r,n}$ congruent to 
$\underline{\alpha}$ 
under $\underline{\pi}_n$, and set $A:=w_\alpha \,\underline{\A}_{r,n}$.
Then 
\[\widehat A = \widehat w \,\underline{\A}_{r,n}\]
where $\widehat w \in W^f$
is given by 
$
\widehat w = \underline{w}(\widehat\mu + \rho_r,n)\,w_{0,\xi},
$
and $w_{0,\xi}$ is the longest element of the stabilizer of
$\xi$ in $\ASG_r$.
\end{lemma}
\Proof By definition of $\mu^{(0)}$ and $\mu^{(1)}$ (see Section~\ref{SECT2}),
the translated alcove $B:=A-n\underline{\mu^{(1)}}$ lies in $n\Pi$,
where 
\[
n\Pi =\{ \underline{\lambda} \in \underline{P}_r \ | \ 
0\le \lambda_i-\lambda_{i+1} \le n,\ i=1,\ldots , r\}
\]
is the fundamental box for the level $n$ action of $\ASG_r$.
Hence, by definition of $\widehat A$ (see \cite{So1}), we
can write 
$\widehat A = w_0(B) + 2n\underline{\rho}_r + n\underline{\mu^{(1)}}$,
which shows that $\widehat A$ contains $\underline{\widehat{\mu}+\rho_r}$.
Thus, if the stabilizer of $\xi$ is trivial, we are done. 
Otherwise, since $w_\alpha$ is of minimal length
among the elements $w\in\ASG_r$ such that $w \xi = \underline{\alpha}$,
we see that $A$ is the lowest alcove adjacent to $\underline{\alpha}$,
and therefore $\widehat A$ is the highest alcove 
adjacent to $\underline{\widehat{\mu}+\rho_r}$, which means
that $\widehat A = \widehat w \,\underline{\A}_{r,n}$ with
$\widehat w = \underline{w}(\widehat\mu + \rho_r,n)w_{0,\xi}$.
\cqfd

\begin{proposition}\label{P2}
Let $\lambda ,\,\mu$ be partitions of $m$ with at most $r$ parts.
Put $\beta = \lambda+\rho_r$, $\alpha=\mu+\rho_r$.
Then, with the notation of Theorem~\ref{TH2},
\[
r_{\beta,\alpha}(q) = 
q^{\ell-\ell_\mu}\,d_{\widetilde\lambda,\widehat\mu}(q^{-1})\,.
\]
\end{proposition}
\Proof
Let $\xi=w_\alpha^{-1}\underline{\alpha}$.
For $s\in\SG_r$ we have
\[
Q_{sw_\beta\,,\,w_0w_\alpha}(q)
=(-q)^{-\ell(w_{0,\xi})}\,
Q_{sw_\beta w_{0,\xi}\,,\,w_0w_\alpha}(q)\,
\]
(see \cite{So1}, proof of Proposition~3.7).
Hence, we can rewrite Proposition~\ref{PROP5} as 
\[
r_{\beta,\alpha}(q) = (-q)^{-\ell(w_{0,\xi})}\, 
m^{w_\beta w_{0,\xi}\, ,\, w_\alpha}(q)\,.
\]
By \cite{So1}, Theorem~5.1, and Lemma~\ref{LEMM} above, we get
\[
m^{w_\beta w_{0,\xi}\, ,\, w_\alpha}(q)
=
q^{\ell(w_0)}\,n_{w_\beta w_{0,\xi}\,,\,
\underline{w}(\widehat\mu+\rho_r,n)w_{0,\xi}}
(q^{-1})\,,
\]
where for $x,y\in W^f$,
\[
n_{x,y}(q) = \sum_{s\in\SG_r} (-q)^{\ell(s)}\,P_{sx,y}(q)\,,
\]
(\cite{So1}, Prop. 3.4).
On the other hand, the expression of $d_{\lambda,\mu}(q)$ in terms
of Kazhdan-Lusztig polynomials given by Varagnolo-Vasserot may be written as 
\[
d_{\lambda,\mu}(q) = \sum_{s\in\SG_r} (-q)^{\ell(s)}
P_{sw_\beta w_{0,\xi}\,,\,w_\alpha w_{0,\xi}}(q)
\]
(see \cite{LT2}, Eq.~(93)).
Hence 
\begin{equation}\label{TILT}
d_{\lambda,\mu}(q) = n_{w_\beta w_{0,\xi} \, , \, 
w_\alpha w_{0,\xi}}(q)\,,
\end{equation}
and using the fact that $\underline{\lambda} = \underline{\widetilde{\lambda}}$, we get
\[
r_{\beta,\alpha}(q) = (-q)^{\ell(w_0)-\ell(w_{0,\xi})}\,
n_{w_\beta w_{0,\xi}\,,\,\underline{w}(\widehat\mu+\rho_r,n)w_{0,\xi}}
(q^{-1})
= (-q)^{\ell(w_0)-\ell(w_{0,\xi})}\,
d_{\widetilde{\lambda},\widehat{\mu}}(q^{-1}) \,.
\]
\cqfd

Combining Proposition~\ref{P1} and Proposition~\ref{P2},
we have proved Theorem~\ref{TH2}.

\section{Lusztig's and Soergel's formulas}
\label{SECT5}

In \cite{LT} (see also \cite{LT2}), the polynomials $d_{\lambda,\mu}(q)$
and $e_{\lambda,\mu}(q)$ were introduced as the coefficients of
two canonical bases $\{\G^+_\lambda\}$ and $\{\G^-_\lambda\}$
of the Fock space $\F$, and it was conjectured that
\begin{eqnarray}
d_{\lambda,\mu}(1) &=& [W(\lambda') : L(\mu')] = d_{\lambda,\mu}\,, \label{C1} \\
e_{\lambda,\mu}(-1) &=& [L(\mu) : W(\lambda)]\,. \label{C2}
\end{eqnarray}
(The notation of \cite{LT2} slightly differs from that of \cite{LT}, and
here we follow \cite{LT2}.)
Let $S(\lambda)$ denote the Specht module for $H_m(v)$ corresponding
to $\lambda$.
If $\mu$ is $n$-regular, $S(\mu)$ has a simple head that we denote
by $D(\mu)$, and it is known \cite{DJ} that
\[
[S(\lambda) : D(\mu)] = d_{\lambda,\mu} \,.
\]
The subset $\{\G^+_\lambda \ | \ \mbox{$\mu$ is $n$-regular}\}$
is nothing but Kashiwara's global crystal basis \cite{Ka2} of
$U_q^-(\slchap_n)\,|\emptyset \> \cong L(\Lambda_0)$, and it had 
been previously conjectured in \cite{LLT} that for $\mu$ $n$-regular
\begin{equation}\label{C0}
d_{\lambda,\mu}(1) = [S(\lambda) : D(\mu)]\,.
\end{equation}
This conjecture was proved by Ariki \cite{Ar} (see \cite{Ge}
for a detailed review of this work).

Using Theorem~\ref{TH1} and Theorem~\ref{TH2} we obtain immediately
that (\ref{C0}) implies (\ref{C1}), that is, Ariki's theorem yields
also a proof of the main conjecture of \cite{LT}.

On the other hand, we know that the matrix
$[e_{\lambda',\,\mu'}(-q)]$ is the inverse of 
$[d_{\lambda,\mu}(q)]$ (see \cite{LT} and \cite{LT2}, Theorem~7.14).
Hence (\ref{C1}) implies (\ref{C2}).

It was proved by Varagnolo-Vasserot \cite{VV} that the $e_{\lambda,\mu}(q)$
are exactly the parabolic Kazhdan-Lusztig polynomials occuring
in Lusztig's formula for the expression of the character of $L(\mu)$
in terms of those of the $W(\lambda)$.
Therefore, we see that Lusztig's formula can be derived from
Ariki's theorem.

Finally, using again the following result of Du-Parshall-Scott
and Donkin (see Section~\ref{SECT3})
\[
[T(\mu) : W(\lambda) ] = [W(\lambda') : L(\mu')]
= d_{\lambda,\mu}\,,
\]
and the expression (\ref{TILT}) of $d_{\lambda,\mu}(q)$ as a parabolic
Kazhdan-Lusztig polynomial obtained by Varagnolo-Vasserot 
(see also Goodman-Wenzl \cite{GW}), we recover Soergel's
formula for the character of the tilting module $T(\mu)$.


\subsection*{Acknowledgements}
I thank R. Rouquier, M. Geck for stimulating discussions,
and H. H. Andersen for pointing out Ref.~\cite{AnP}.
This work was started during my stay at R.I.M.S. Kyoto University
(07/98-11/98) and my visit at the University of Tottori (09/98).
I wish to thank the organizers of the R.I.M.S. Project 1998
``Combinatorial methods in representation theory'',
K. Koike, M. Kashiwara, S. Okada,  H.-F. Yamada, I. Terada
for their invitation, M. Ishikawa for his warm hospitality in Tottori,
and the Hanabi people for their kindness.



\bigskip\bigskip


\begin{thebibliography}{ABC} \scriptsize
%
\bibitem{An} {\sc H. H. Andersen}, {\it Tensor products
of quantized tilting modules}, Commun. Math. Phys. {\bf 149}
(1992), 149-159.
%
\bibitem{AnP} {\sc H. H. Andersen, J. Paradowski}, 
{\it Fusion categories arising from semisimple Lie algebras},
Commun. Math. Phys. {\bf 169} (1995), 563-588.
%
\bibitem{APW2}{\sc H. H. Andersen, P. Polo, K. Wen},
{\it Injective modules for quantum algebras}, 
Am. J. Math. {\bf 114} (1992), 571-604.
%
\bibitem{Ar}{\sc S. Ariki},
{\it On the decomposition numbers of the Hecke algebra of $G(m,1,n)$},
J. Math. Kyoto Univ. {\bf 36} (1996), 789-808.
%
\bibitem{DJ} {\sc R. Dipper, G. D. James}, {\it The $q$-Schur algebra},
Proc. London Math. Soc. {\bf 59} (1989), 23-50.
%
\bibitem{Do} {\sc S. Donkin}, {\it The $q$-Schur algebra},
L.M.S. Lect. Notes {\bf 253}, Cambridge Univ. Press, 1998.
%
\bibitem{Du} {\sc J. Du}, {\it A note on quantized Weyl reciprocity
at roots of unity}, Algebra Colloq. {\bf 2} (1995), 363-372.
%
\bibitem{DPS} {\sc J. Du, B. Parshall, L. Scott},
{\it Quantum Weyl reciprocity and tilting modules},
Commun. Math. Phys. {\bf 195} (1998), 321-352.
%
\bibitem{Erd1} {\sc K. Erdmann}, {\it Decomposition numbers for
symmetric groups and composition factors of Weyl modules},
J. Algebra {\bf 180} (1996), 316-320.
%
\bibitem{Ge} {\sc M. Geck}, {\it Representations of Hecke algebras
at roots of unity}, S\'eminaire Bourbaki, 1997-98, expos\'e 836.
%
\bibitem{GW} {\sc F. Goodman, H. Wenzl}, {\it Crystal bases of
quantum affine algebras and affine Kazhdan-Lusztig polynomials},
math.QA/9807014.
%
\bibitem{JK} {\sc G. D. James, A. Kerber}, {\it The representation theory of
the symmetric group}, Addison-Wesley, 1981.
%
\bibitem{Ja1} {\sc G. D. James}, {\it The irreducible representations
of the finite general linear groups},
Proc. London Math. Soc. {\bf 52} (1986), 236-268.
%
\bibitem{Jam} {\sc G. D. James}, {\it The decomposition matrices
of $GL_n(q)$ for $n\le 10$},
Proc. London Math. Soc., {\bf 60} (1990), 225-265.
%
\bibitem{Ja}{\sc J. C. Jantzen}, {\it Representations of
algebraic  groups},
Academic Press, 1987.
%
%
\bibitem{Ka2}{\sc M. Kashiwara}, {\it On crystal bases of the $q$-analogue
of universal enveloping algebras},  Duke Math. J. {\bf 63} (1991), 465-516.
%
\bibitem{KT}{\sc M. Kashiwara, T. Tanisaki}, {\it Kazhdan-Lusztig conjecture
for affine Lie algebras with negative level},
Duke Math. J. {\bf 77} (1995), 21-62.
%
%
\bibitem{KL}{\sc D. Kazhdan, G. Lusztig}, {\it Tensor structures arising
from affine Lie algebras I-II, III-IV},
J. Am. Math. Soc. {\bf 6\ 7 }, (1994) 905-1011; (1994) 335-453.
%
\bibitem{LLT}{\sc A.~Lascoux, B.~Leclerc, J.-Y.~Thibon},
{\it Hecke algebras at roots of unity and crystal bases of
quantum affine algebras},
Commun. Math. Phys. {\bf 181} (1996), 205-263.
%
\bibitem{LT} {\sc B. Leclerc, J.-Y. Thibon}, {\it Canonical bases
of $q$-deformed Fock spaces}, Int. Math. Res. Notices, {\bf 9} (1996),
447-456.
%
\bibitem{LT2}{\sc B.~Leclerc, J.-Y.~Thibon},
{\it Littlewood-Richardson coefficients and Kazhdan-Lusztig polynomials},
math.QA/9809122.
%
\bibitem{Lu5} {\sc G. Lusztig}, {\it On quantum groups},
J. Algebra, {\bf 131} (1990), 466-475.
%
%
\bibitem{OS}{\sc O. Schiffmann}, {\it Personal communication}.
%
\bibitem{So1} {\sc W. Soergel}, {\it Kazhdan-Lusztig-Polynome und
eine Kombinatorik f\"ur
Kipp-Moduln}, Represent. Theory 1 (1997), 37-68 (english 83-114).
%
\bibitem{So2} {\sc W. Soergel},  {\it Charakterformeln f\"ur
Kipp-Moduln \"uber Kac-Moody-Algebren},
Represent. Theory 1 (1997), 115-132.
%
\bibitem{VV} {\sc M. Varagnolo, E. Vasserot}, {\it On the decomposition 
matrices of the quantized Schur algebras}, math.QA/9803023.
%

\end{thebibliography}
\end{document}